\documentclass[]{amsart}

\usepackage{amsfonts,amssymb,amsmath}
\usepackage[pdftex]{graphics}
\usepackage{pdfsync}

\begin{document}

\hyphenation{
sub-se-quen-ces
}

\title{On the notion  of a basis of a finite dimensional vector space}

\author{Alexander Gamkrelidze}
\address{Department of Computer Science\\
I. Javakhishvili Tbilisi State University\\
Georgia
}
\email{alexander.gamkrelidze@tsu.ge}
\author{Grigori Giorgadze}
\address{Department of Mathematics\\
I. Javakhishvili Tbilisi State University\\
Georgia}
\email{gia.giorgadze@tsu.ge}

\date{}

\maketitle

 \section*{Abstract}
In this Note, we show that the notion of a basis of a finite-dimensional vector space could be introduced by an argument much weaker than Gauss' reduction method. Our aim is to give a short proof of a simply formulated lemma, which in fact is equivalent to the theorem on frame extension, using only a simple notion of the kernel of a linear mapping, without any reference to special results, and derive the notions of basis and dimension in a quite intuitive and logically appropriate way, as well as obtain their basic properties, including a lucid proof of Steinitz's theorem.

 \section{Introduction}
 To indroduce the notion of a basis of a finite dimensional vector space the 
following Proposition or its equivalent should be proved:

\emph{Given two $n-$frames
(sequence of linearly independent vectors) 
 $\|e_j\|_n, \|f_j\|_n$  in a vector space, the vector
$f_j$ contained in the linear span of $\|e_j\|_n$, $f_j\in [\|e_j\|_n], j=1,\ldots,n$,} 
\emph{there exists a uniquely defined}  \emph{invertible  $n\times n$-matrix of scalars $A$, satisfying  the equations}
\[
  \|f_j\|_n=\|e_j\|_nA \iff \|e_j\|_n=\|f_j\|_n \,A^{-1}. 
\]

\noindent
Existence and uniqueness of the matrix $A$ is evident, but proving its inverti\-bility is always a delicate pedagogical problem for every author of an introduc\-tory text on Linear Algebra, since it should be given right at the beginning of the exposition, after the vector space axioms are listed. 

\noindent
A standard argument of proving the Proposition is based
 (explicitly or implicitly)
  on Gauss' reduction process  of transforming a square matrix   to a 
triangular
  form, often combined with the Steinitz's theorem on   frame extension. 

\noindent
In this Note we suggest, instead of solving the above equations, to solve the system of inclusions
\[
e_i\in[\|f_j\|_n], i=1,\ldots, n,
\]
which is equivalent  under given conditions to  the above equations.    The proof is short and based on an argument much weaker 
than
the reduction process - additionally to the list of axioms we only need the notion of the kernel of a linear mapping, without any special constructions performed on the data.

\

\noindent
The authors thank Prof. Revaz V. Gamkrelidze for drawing our attention to the problem and stimulating discussions during the preparation of this manuscript.

\

 \section{Preliminary remarks}
 The logical procedure of introducing   a basis of a finite dimensional vector space
is completely equivalent to that of inverting a nondegenerate square matrix of scalars and is reduced to proving the following 
proposition
or its equivalent, after which the path to the notions of a basis and dimension 
are almost uniquely determined and short.
\vskip0,05in
\noindent
\underline  {Proposition.} \emph{If two n-frames of a vector space (not necessarily finite dimensional), $\|e_j\|_n$ and $\|f_j\|_n$,  are given, where the vectors $f_j, \ j=1,\ldots, n$, are con\-tained in the linear span   
$[\|e_j \|_n]$, then there exists an invertible $n\times n$   matrix of scalars $A$ such that }
\begin{equation}
\label{1}
\|f_j\|_n=\|e_j\|_n A   \iff  \|e_j\|_n =\|f_j\|_nA^{-1}.    
 \end{equation}

\noindent
 The existence and uniqueness of   $A$ are evident since the sequences $\|e_j\|_n$, $\|f_j\|_n$ are linearly independent, but proving its invertibility is always a deli\-cate pedagogical problem for every author of an introductory text on Linear Algebra, since the proof should be given right at the beginning of the exposi\-tion, after the vector space axioms are listed, 
 very often  {even} before linear mappings are introduced.
 
 \noindent
 A standard way of solving the problem is  to apply  to $A$ (implicitly or explici\-tly) Gauss' method of reduction of a square matrix to a triangular form, often combi\-ned with Steinitz's theorem on a frame extension. In this respect, it is instructive to compare   expositions of the  corresponding material in four texts 
 \cite{Algebra1}, \cite{Algebra2}, \cite{Algebra3} and \cite{Algebra4}, the last three published after \cite{Algebra1} with the time lag of more than 50 years.

\noindent
In this 
{
 {note,}}    
instead of solving   equations (1),   we prove the system of inclu\-sions
\begin{equation}
e_i\in [\|f_j\|_n],\ i=1,\ldots,n,
\end{equation}
  which is equivalent to (1) since $\|e_j\|_n, \|f_j\|_n$ are linearly independent.
 The proof is much shorter and   based on an argument much weaker than the reduction process. 
  {Additionally}  
to the list of axioms we only need the notion of the kernel of a linear mapping, without any special constructions performed on the data.

\noindent
{

In section 3, inclusions (2) are presented and discussed as a short and simply formulated Basic Lemma, from which, in sections 4 - 6, the notions of a basis and dimension of a finite dimensional vector space 
along with their basic properties (including a lucid proof of the frame extension theorem) 
are derived in a quite intuitive and logically proper way. Finally, in Section 7, the Basic Lemma is proved. In the next Section 2 we list several necessary initial notions used in this note. 	
}

\noindent
To conclude the section, we should remark that in most texts on introductory Linear Algebra, linear mappings are introduced  too late, contrary to a   well known motto, according to which ``morphisms in a category are at least as useful as the objects are''.   We think that   linear mappings and   their basic properties should be exposed in introductory texts   right after the   vector space axioms are listed   and   properly discussed, and consider this Note as a proper example supporting our opinion.

 \section{A short list of necessary initial notions} 
 Together with the standard general set-theoretic notions related to a map\-ping $L$,
 \[
 dom\, L=\mathbb E,\  codom\, L=\mathbb F,\ im\,L= L(\mathbb E),
 \]
  it is useful to have at the disposal from the very beginning  
    {  specifical} 
   linear notions of a (linear) subspace $\mathbb F\subset\mathbb E$ and the factor $\mathbb E\,/\mathbb F$, in particular, the notions of the kernel $ker\,L$ and   factor $\mathbb E/ ker\,L$. 
 
 \noindent 
For an adequate definition of the basis 
 and its appropriate discus\-sion,  we should have a certain freedom in handling finite sequences of vectors in $\mathbb E$.  An arbitrary sequence of length $n$ of vectors in $\mathbb E$   is presented as an $n$-row matrix,
   \[
   ||x_1,\ldots,x_n||=\|x_j\|_n=\|x_j\|,\ x_j\in\mathbb E,\ j=1,\ldots,n.
   \]
 Every sequence of vectors could be \emph{extended } by ascribing to it new vectors from the space. It is useful to remember that a sequence of length $n$ of vectors in $\mathbb E$ is not a subset of $\mathbb E$, but rather a function to $\mathbb E$ on the ordered set of first $n$ naturals, where  $im \|x_j\|\subset\mathbb E$. If $im\|x_j\|$ is a subset of a subspace $\mathbb F\subset\mathbb E$ we say that the sequence $\|x_j\|$ belongs to  $\mathbb F$ and write $\|x_j\|\prec\mathbb F$.
 
 The initial list of basic items should  also contain the following notions: linear combination of vectors of a sequence, linearly independent and depen\-dent sequences of vectors of $\mathbb E$, $n$ - frames
  --- linearly independent sequences of length $n$, or \emph{rank} $n$, 
  \[
    \|e_j\|=\|e_j\|_n,\   rank\, \|e_j\|_n=n.
  \]
  We also need the following notions:
the rank of an arbitrary sequence,\break 
$rank\,\|x_j\|,\ \|x_j\| \prec\mathbb E$,  defined as the maximal   rank  of frames contained (as subsequences) in $\|x_j\|$, the linear span $[\|x_j\|]\subset\mathbb E$ of an arbitrary sequence of vectors $\|x_j\|\prec\mathbb E$, linear span of an arbitrary subset $A\subset\mathbb E$ --- minimal subspaces in $\mathbb E$ (with respect to the set-theoretic inclusion) containing, respec\-tively, subsets $im\|x_j\|$ and $A$.

 \section{The basic lemma. Formulation and a preliminary discussion} 
   \noindent A frame $\|e_j\|\prec \mathbb F\subset\mathbb E$ is \emph{maximal in a subspace} $\mathbb F$ if the rank of an arbitrary sequence $\|x_j\|\prec \mathbb F$ is bounded by the rank of $\|e_j\|$,
 \[
 \|x_j\|\prec\mathbb F\Longrightarrow rank\, \|e_j\|\ge rank\, \|x_j\|.
 \]
A frame $\|e_j\|_n\prec\mathbb F$ could be \emph{extended} (is \emph{extendable}) in $\mathbb F$,  if there exists a vector $f\in\mathbb F$ such that the extended sequence $e_1,\ldots,e_n,f$ is again a frame (of rank $n+1$).
Intuitively,  maximality of a frame could be considered as a ``generalized version'' of the ``individual property'' of a frame to be non-extendable. Every maximal frame in $\mathbb F$ is evidently not extendable. The inversion of the assertion is also true --- every non-extendable frame in $\mathbb F$ is maximal in $\mathbb F$ (Steinitz's extension theorem), but the proof is not trivial and in fact  belongs right to the core of the  problem under discussion --- of giving  proper definitions of a basis and dimension of a finite dimensional vector space.

\noindent
 { Now we} 
shall formulate 
 a simple lemma, which easily clears up all inter\-relations between the notions of maximality of a frame and its ability ``to be extendable'', and suggests     natural intuitive definitions of a basis and dimension.

\noindent
 {\bf The Basic Lemma.} \emph {Every frame} $\|e_j\|\prec\mathbb E$ \emph{is maximal in its linear span}
$[\|e_j\|]$ \emph{and extendable in every subspace} $\mathbb F\supset [\|e_j\|] $,  \emph{if the inclusion is strong}.

\noindent
The ability of being ``extendable''  under the given conditions is evident --- it is achieved by ascribing to $\|e_j\|$ 
 an arbitrary vector $f\in\mathbb F, f\notin [\|e_j\|]$,  the maximality of $\|e_j\|$ in the linear span $[\|e_j\|]$ is equivalent to each of the following two assertions.

1. \emph{Every frame $\|f_j\|_n$ of rank $n$    in the linear span $[\|e_j\|_n]$ is not extendable there, or, equivalently, the following system of $n$ inclusions is valid,}
\begin{equation}
\label{INCL}
e_i\in[\|f_j\|_n], \ i=1,\ldots,n.
\end{equation}

2. \emph{The rank of an arbitrary sequence  in $\mathbb E$  consisting of linear combinations of a fixed sequence of $n$ vectors from $\mathbb E$ does not exceed $n$.}

\noindent
A simple inductive proof of the system of inclusions (3), based on the notion of the kernel of a linear mapping, is given (as already stated)    in the final Sec\-tion 6. Before, in next two Sections, we shall derive    
simple and intuitive definitions of    a basis and dimension of a finite dimensional vector space
{from the formulated lemma}, 
as well as give a lucid   proof of   Steinitz's  theorem on the  frame extension.

\section{Definition of a basis and dimension of a finite dimensional vec\-tor space} 
A vector space $\mathbb E$ is \emph{finite dimensional} if it contains a \emph{finite set of generators} --- a finite subset $G\subset\mathbb E$  with the linear span coinciding with $\mathbb E$.   
Since $G$ is finite, it contains a   frame $\|e_j\|_n\prec G$ (of maximal rank in $G$), which spans the whole space $\mathbb E=[\|e_j\|_n]$, hence, according to the lemma, the frame is maximal in $\mathbb E$.
Hence every finite dimensional vector space contains maximal frames --- \emph{the bases of the space}.
Their common rank $n$ is the \emph{dimension of the space}, and every vector $x\in\mathbb E$ 
is uniquely represented as
\[
x=\sum_1^n\lambda^\alpha e_\alpha,\quad \lambda^j\in\Lambda,\ j=1,\ldots,n.
\]
\emph{Thus, every basis of a finite dimensional vector space is an irreducible system of generators of the space.}

Conversely, if a sequence $\|x_j\|_n$ in a vector space $\mathbb E$ is given such that every vector  
$x\in\mathbb E$ is uniquely represented as
\[
x=\sum_1^n\lambda^\alpha x_\alpha,\quad \lambda^j\in\Lambda, 
\]
then, according to the lemma, $\|x_j\|_n$ is a maximal frame of rank $n$, or a basis  in $\mathbb E$.


\section{The Steinitz' theorem on the frame extension} 
 In $k+l$-dimensio\-nal vector space $\mathbb E^{k+l}$ a basis $B=\|e_j\|_{k+l}$ and an arbitrary $k$-frame $\|f_j\|_k$ are given. From general considerations it easily follows by induction that $B$ contains a subsequence $e_{i_1},\ldots,e_{i_r}$ of a \emph{certain length} $r\le k+l$ such that the frame   $\|f_j\|_k$ extended by the subsequence is a frame of length $k+r$,
\[
B'=\|f_1,\ldots,f_k;e_{i_1},\ldots,e_{i_r}\|\prec\mathbb E^{k+l},
\]
containing the linear span  of $B$, which is the whole space $\mathbb E^{k+l}$,  hence
\[
[B']=\mathbb E^{k+l}.
\]
According to the lemma, the obtained frame $B'$ is maximal in $\mathbb E^{k+l}$, i. e. is a basis of the $k+l$-dimensional vector space $\mathbb E^{k+l}$, hence $r=l$, and we have extended the initial $k$-frame $\|f_j\|_k$ by a subsequence of length $r=l$ of a preassigned basis $B$ to a new basis $B'$   of the space $\mathbb E^{k+l}$.

\section{The basic lemma. Proof} 

We shall prove the system of inclusions (2) by induction performed on the rank $n$ of the frame $\|e_j\|_n$. The assertion is evident for $n=1$ since 
{in this case}
the linear span  
coincides with the family of vectors

\[
[\|e_j\|_1]=\left\{\lambda e_1|\, \lambda\in\Lambda\right\}.
\]
Assuming that the assertion is proved for all natural $k\le n-1$, consider the case $k=n$. Introduce $n$ linear mappings
\[
L_i\in Hom\left([\|e_j\|_n], [\|f_1,\ldots,\hat f_i,\ldots,f_n\|]\right),\quad i=1,\ldots,n,
\]
by defining $L_i$ on vectors $e_1,\ldots,e_n$ according to the equations
\[
L_ie_j=f_j,\ i\neq j,\ \ L_ie_i=0,\quad i,j=1,\ldots,n.
\]
The kernel of $L_i$ coincides with the subspace
\[
ker L_i=\{\lambda e_i\bigr|\lambda\in\Lambda\},
\]
and the image  ---  with the linear span
\[
 im\,L_i=[\|f_1,\ldots,\hat f_i,\ldots f_n\|].
\]
The restriction of $L_i$ on the subspace $[\|f_j\|_n]\subset [\|e_j\|_n]$ has a nonzero kernel since otherwise the $n$-frame $\|L_if_j\|, j=1,\ldots,n$ would be embedded in the linear span of an $(n-1)$-frame $\|f_1,\ldots,\hat f_i,\ldots,f_n\|$, which contradicts the inductive assumption. Hence,
\[
\lambda e_i\in ker L_i\Bigr|_{[\|f_j\|_n]} \  \forall\lambda\in\Lambda,\ i=1,\ldots,n,
\]
or
\[
e_i\in[\|f_j\|_n]\ \forall i=1,\ldots,n.
\]
This proves the Lemma.


\begin{thebibliography}{4}


\bibitem{Algebra1}
G. Schreier, E. Sperner.
Einf\"uhrung in die Analytische Geometrie und Algebra,  \linebreak
Bd. 1
B.G. Teubner, Leipzig, 1931




\bibitem{Algebra2}
Ren\`e Deheuvels.
Formes quadratiques et groupes classiques \\
Presses Universitaires de France, 1981


\bibitem{Algebra3}
L\'aszl\'o Babai and 
P\'eter Frankl.
Linear Algebra Methods in Combinatorics,\\
Dept. of Computer Science, The University of Chicago,\\
Preliminary Version, 1992




\bibitem{Algebra4}
Sheldon Axler.
Linear Algebra Done Right\\
Springer Verlag, 2015








\end{thebibliography}
\end{document}